\documentclass[11pt]{article}
\usepackage{amssymb,amsfonts,amsmath,curves}
\usepackage{epsfig}
\usepackage{graphicx}
\usepackage{color}
\usepackage{anysize}
\usepackage{hyperref}
\usepackage{url}

\marginsize{2cm}{2cm}{1cm}{1cm}

\begin{document}

\makeatletter\@addtoreset{equation}{section}
\makeatother\def\theequation{\thesection.\arabic{equation}}

\baselineskip=14pt

\newtheorem{defin}{Definition}[section]
\newtheorem{teo}{Theorem}[section]
\newtheorem{ml}{Main Lemma}
\newtheorem{con}{Conjecture}
\newtheorem{cond}{Condition}
\newtheorem{conj}{Conjecture}
\newtheorem{prop}[teo]{Proposition}
\newtheorem{lem}{Lemma}[section]
\newtheorem{rmk}[teo]{Remark}
\newtheorem{cor}{Corollary}[section]

\newcommand{\be}{\begin{equation}}
\newcommand{\ee}{\end{equation}}
\newcommand{\bp}{\begin{prop}}
\newcommand{\ep}{\end{prop}}
\newcommand{\bt}{\begin{teo}}
\newcommand{\bcor}{\begin{cor}}
\newcommand{\ecor}{\end{cor}}
\newcommand{\bcon}{\begin{con}}
\newcommand{\econ}{\end{con}}
\newcommand{\bcond}{\begin{cond}}
\newcommand{\econd}{\end{cond}}
\newcommand{\bconj}{\begin{conj}}
\newcommand{\econj}{\end{conj}}
\newcommand{\et}{\end{teo}}
\newcommand{\brm}{\begin{rmk}}
\newcommand{\erm}{\end{rmk}}
\newcommand{\bl}{\begin{lem}}
\newcommand{\el}{\end{lem}}
\newcommand{\ben}{\begin{enumerate}}
\newcommand{\een}{\end{enumerate}}
\newcommand{\bei}{\begin{itemize}}
\newcommand{\eei}{\end{itemize}}
\newcommand{\bdf}{\begin{defin}}
\newcommand{\edf}{\end{defin}}
\newcommand{\bpr}{\begin{proof}}
\newcommand{\epr}{\end{proof}}

\newenvironment{proof}{\noindent {\em Proof}.\,\,}{\hspace*{\fill}$\halmos$\medskip}

\newcommand{\halmos}{\rule{1ex}{1.4ex}}
\def \qed {{\hspace*{\fill}$\halmos$\medskip}}

\newcommand{\fr}{\frac}
\newcommand{\Z}{{\mathbb Z}}
\newcommand{\R}{{\mathbb R}}
\newcommand{\E}{{\mathbb E}}
\newcommand{\C}{{\mathbb C}}
\renewcommand{\P}{{\mathbb P}}
\newcommand{\N}{{\mathbb N}}
\newcommand{\Q}{{\mathbb Q}}
\newcommand{\var}{{\mathbb V}}
\renewcommand{\S}{{\cal S}}
\newcommand{\T}{{\cal T}}
\newcommand{\W}{{\cal W}}
\newcommand{\X}{{\cal X}}
\newcommand{\Y}{{\cal Y}}
\newcommand{\h}{{\cal H}}
\newcommand{\Fi}{{\cal F}}
\newcommand{\Mi}{{\cal M}}
\newcommand{\Ni}{{\cal N}}

\renewcommand{\a}{\alpha}
\renewcommand{\b}{\beta}
\newcommand{\g}{\gamma}
\newcommand{\G}{\Gamma}
\renewcommand{\L}{\Lambda}
\renewcommand{\l}{\lambda}
\renewcommand{\d}{\delta}
\newcommand{\D}{\Delta}
\newcommand{\e}{\epsilon}
\newcommand{\eps}{\epsilon}
\newcommand{\s}{\sigma}
\newcommand{\B}{{\cal B}}
\renewcommand{\o}{\omega}

\newcommand{\nn}{\nonumber}
\renewcommand{\=}{&=&}
\renewcommand{\>}{&>&}
\newcommand{\<}{&<&}
\renewcommand{\le}{\leq}
\newcommand{\+}{&+&}

\newcommand{\pa}{\partial}
\newcommand{\ffrac}[2]{{\textstyle\frac{{#1}}{{#2}}}}
\newcommand{\dif}[1]{\ffrac{\partial}{\partial{#1}}}
\newcommand{\diff}[1]{\ffrac{\partial^2}{{\partial{#1}}^2}}
\newcommand{\difif}[2]{\ffrac{\partial^2}{\partial{#1}\partial{#2}}}

\title{A Monotonicity Result for the Range of a Perturbed Random Walk}
\author{Lung-Chi Chen$^{\,1}$ \and Rongfeng Sun$^{\,2}$}
\date{Nov 08, 2012}
\maketitle

\footnotetext[1]{Department of Mathematics, Fu-Jen Catholic University, 510 Chung Cheng Road, Hsinchuang , Taipei County 24205, Taiwan.
E-mail: lcchen@math.fju.edu.tw}

\footnotetext[2]{Department of Mathematics, National University of Singapore, 10 Lower Kent Ridge Road, 119076 Singapore. Email: matsr@nus.edu.sg}

\begin{abstract}
We consider a discrete time simple symmetric random walk on $\Z^d$, $d\geq 1$, where the path of the walk is perturbed by inserting deterministic jumps. We show that for any time $n\in\N$ and any deterministic jumps that we insert, the expected number of sites visited by the perturbed random walk up to time $n$ is always larger than or equal to that for the unperturbed walk. This intriguing problem arises from the study of a particle among a Poisson system of moving traps with sub-diffusive trap motion. In particular, our result implies a variant of the Pascal principle, which asserts that among all deterministic trajectories the particle can follow, the constant trajectory maximizes the particle's survival probability up to any time $t>0$.

\vspace{12pt}

\end{abstract}

\section{Introduction}
Let $Z:=(Z_n)_{n\geq 0}$ be a discrete time random walk on $\Z^d$ with increment distribution $p(\cdot)$, which we assume to be symmetric, i.e., $p(x)=p(-x)$ for all $x\in\Z^d$. We ask the following question: what will happen to the expected number of sites visited by $Z$ up to time $n\in\N$, if deterministic jumps are inserted in $Z$ (at half integer times)? A natural conjecture is that the expected number of visited sites can only increase when jumps are inserted. However, verifying this conjecture even for the simple symmetric random walk on $\Z^d$ for any $d\geq 1$ turns out to be far from trivial, which is the main result of this paper. In dimension $1$, we are able to go beyond the simple symmetric random walk to a more general class of symmetric random walks. The case for general symmetric random walks on $\Z^d$, $d\geq 1$, remains open.

We now formulate our main result. To simplify notation, we will change time units such that the random walk jumps at even integer times, while deterministic jumps are inserted at odd integer times.

\bt\label{T:range}{\bf [Range of a random walk under insertion perturbation]}
Let $(\bar Z_n)_{n\geq 0}$ be a random path on $\Z^d$ such that $(\bar Z_{2k})_{k\geq 0}$ is a random walk with increment distribution $p(\cdot)$, and $\bar Z_{2k}=\bar Z_{2k+1}$ for all $k\geq 0$. Let $R_n(\bar Z):=\{\bar Z_i: 0\leq i\leq n\}$ denote the range of $\bar Z$ up to time $n$, and let $|R_n(\bar Z)|$ denote its cardinality. Assume that $p(\cdot)$ falls into one of the following three classes:
\begin{itemize}
\item[\rm (i)] $p(\cdot)$ is any symmetric distribution on $\Z$ with $p(k)\geq p(k+1)$ for all $k\geq 1$, and $p(0)\geq p(3)$;

\item[\rm (ii)] $p(\cdot)$ is the increment distribution of a simple symmetric random walk on $\Z^d$, for any $d\geq 1$;

\item[\rm (iii)] $p(\cdot)$ is any symmetric distribution on $\Z^d$, $d\geq 1$, with $p(0)\geq \frac{1}{2}$.
\end{itemize}
Then for any path $(f_n)_{n\geq 0}$ on $\Z^d$ with $f_{2k-1}=f_{2k}$ for all $k\in\N$, we have
\be\label{range}
\E[|R_{n}(\bar Z)|] \leq \E[|R_{n}(\bar Z+f)|] \qquad \mbox{for all } n\in\N.
\ee
\et
Theorem~\ref{T:range} under condition (iii) in fact follows from an earlier result of Moreau et al in~\cite{MOBC03, MOBC04} ,  which we will explain in more detail later. We include this result here for completeness.

Our original interest lies in the study of a particle among a Poisson field of mobile traps on $\Z^d$, from which the monotonicity question on the range of a perturbed symmetric random walk arises. From a probabilistic point of view, the latter question is very natural and intriguing on its own, which is why we choose it to be the focus of this paper. Here is the trapping problem we were originally interested in. At time $0$, there is $N_y$ number of traps at each $y\in\Z^d$, where $\{N_y\}_{y\in\Z^d}$ are distributed as i.i.d.\ Poisson random variables with mean $1$. Each trap then moves independently as a random walk on $\Z^d$ with i.i.d.\ holding times, with increment distribution $p(\cdot)$ on $\Z^d$ for the jumps and holding time distribution $\mu(\cdot)$ on $(0,\infty)$ for the time between successive jumps.
For each $y\in\Z^d$ and $1\leq j\leq N_y$, let $Y^y_j:=(Y^y_j(t))_{t\geq 0}$ denote the time-evolution of the $j$-th trap starting from $y$ at time $0$. Then at each time $t\geq 0$, the trap configuration is determined by
\be\label{xitx}
\xi(t,x) := \sum_{y\in \Z^d, 1\leq j\leq N_y} 1_{\{Y^y_j(t)=x\}}, \qquad x\in \Z^d.
\ee
We will denote probability and expectation for $\xi$ by $\P^\xi$ and $\E^\xi$ respectively. A particle $X:=(X(t))_{t\geq 0}$ moving on $\Z^d$ is then killed at the first time
\be\label{tauXxi}
\tau_{X,\xi}:=\inf\{t\geq 0: \xi(t, X(t))\geq 1\},
\ee
when the particle first meets a trap. The particle motion $X$ may be either deterministic or random. We are interested in the probability $S_t$ that $X$ survives up to time $t$, when the randomness in $X$ and the trap field $\{Y^y_j\}_{y\in\Z^d, 1\leq j\leq N_y}$ have been averaged out.

For more background on the trapping problem described above, see e.g.\ \cite{MOBC04, YOLBK08, BAY09, S98, DGRS12} and the references therein. When the traps are mobile, it is in general difficult to obtain good upper bounds on $S_t$. One approach developed in the physics literature is the so-called {\em Pascal principle}, which asserts that among all deterministic trajectories the particle $X$ may follow, the constant trajectory maximizes the survival probability $S_t$. A discrete time version of the Pascal principle was established rigorously by Moreau et al in~\cite{MOBC03, MOBC04} for traps which follow independent random walks with increment distribution $p(\cdot)$ that satisfy Theorem~\ref{T:range}~(iii). Their result can be equivalently formulated as the statement that (see~\cite[Sec.~2.4]{DGRS12}): the expected number of sites visited by a random walk $Z$ up to time $n$ can only increase if we replace $(Z_i)_{i\geq 0}$ by $(Z_i+f_i)_{i\geq 0}$ for any deterministic function $f:\N_0\to\Z^d$. It is easy to see that this monotonicity result for the range of $Z$ under such an additive perturbation implies the monotonicity result for the range of $Z$ under insertion perturbation formulated in Theorem~\ref{T:range}, which justifies Theorem~\ref{T:range}~(iii). By discrete time approximation, the Pascal principle established by Moreau et al can then be used to deduce the Pascal principle for the trapping problem where the holding time distribution $\mu$ is {\em exponential} and $p(\cdot)$ is any symmetric distribution on $\Z^d$ (see~\cite[Sec.~2.4]{DGRS12}).

Recently, the trapping problem with sub-diffusive trap motion has been studied in the physics literature~\cite{YOLBK08, BAY09}, where the Pascal principle was assumed to hold and then used to give bounds on the survival probability $S_t$, the decay rate of which
was then found to be different from the case with diffusive trap motions. The sub-diffusive trap motions were modeled by random walks with heavy-tailed holding time distributions. However, when the holding time distribution $\mu$ is not exponential, one can no longer deduce the Pascal principle from the result of Moreau et al. This motivates us to give a rigorous proof of the Pascal principle for the trapping problem with a general continuous holding time distribution. Our investigation led to the monotonicity question on the range of a symmetric random walk under insertion perturbation.

We formulate below the precise Pascal principle we obtain for the trapping problem, which is effectively a corollary of Theorem~\ref{T:range}.

\bt\label{T:Pascal} {\bf [Pascal Principle for continuous time trapping]} Let $\xi$ be defined as in (\ref{xitx}), where the traps' increment distribution $p(\cdot)$ satisfies one of the conditions in Theorem~\ref{T:range}~(i)--(iii), and the holding time distribution $\mu(\cdot)$ is continuous. Then for any $t\geq 0$ and any $X(\cdot) : [0,\infty)\to\Z^d$ with locally finitely many jumps, we have
\be\label{Pascal}
S_t(X):=\P^\xi(\tau_{X,\xi} >t) \ \leq \ S_t(0):=\P^\xi(\tau_{0,\xi} >t),
\ee
where $\tau_{X,\xi}$ is defined in (\ref{tauXxi}), with $\tau_{0,\xi}$ for the case $X(\cdot)\equiv 0$.
\et

The rest of the paper is organized as follows. In Section~\ref{S:reduction}, we will show how Theorem~\ref{T:Pascal} follows from Theorem~\ref{T:range}, which can then be further reduced to a discrete time trapping problem. Theorem~\ref{T:range} is then proved in Section~\ref{S:Z} for a class of symmetric random walks on $\Z$, and proved in Section~\ref{S:SRW} for the simple symmetric random walk on $\Z^d$ for all $d\geq 1$.

\section{Reduction to Discrete Time Trapping}\label{S:reduction}
We first explain how does Theorem~\ref{T:Pascal} follow from Theorem~\ref{T:range}. Let us fix a realization of the particle motion $X(\cdot):[0,\infty)\to\Z^d$ with locally finitely many jumps, as in Theorem~\ref{T:Pascal}. By integrating out the Poisson field $\xi$, we have
\be
S_t(X) = \P^\xi(\tau_{X,\xi}>t) =\prod_{y\in\Z^d} \exp\{-1+  \P^Y_y(\tau_X>t)\} = \exp\Big\{-\sum_{y\in\Z^d} \P^Y_y(\tau_X\leq t)\Big\},
\ee
where $\P^Y_y$ denotes probability for a trap $Y$ starting at $y\in\Z^d$ at time $0$, and
\be
\tau_X:=\tau_X(Y):=\tau_0(Y-X):=\inf\{t\geq 0: Y(t)-X(t)=0\}.
\ee
Therefore (\ref{Pascal}) reduces to
\be\label{target}
\sum_{y\in\Z^d} \P^Y_y(\tau_X\leq t) \geq \sum_{y\in\Z^d} \P^Y_y(\tau_0\leq t).
\ee
By translation invariance, this can be rewritten in terms of a single trap starting from the origin:
$$
\sum_{y\in\Z^d} \P^Y_0(\tau_{-y}(Y-X)\leq t)  \geq \sum_{y\in\Z^d} \P^Y_0(\tau_{-y}(Y)\leq t),
$$
which is equivalent to
\be\label{rangecomp1}
\E^Y_0[ |R_t(Y-X)|] \geq \E^Y_0[ |R_t(Y)|],
\ee
where $R_t(Y-X):=\{Y(s)-X(s): 0\leq s\leq t\}$ is the range of $Y-X$ up to time $t$. Since $Y$ is a random walk with i.i.d.\ holding times, we can condition on the times at which $Y$ jumps. Note that the jumps of $Y$ and $X$ almost surely do not occur at the same time because the holding time distribution $\mu(\cdot)$ is continuous, therefore $Y-X$ has the effect of inserting the jumps of $-X$ into $Y$. Thus (\ref{rangecomp1}) would follow once Theorem~\ref{T:range} is established.

\begin{figure}[tp] 
\begin{center}
\includegraphics[width=11cm]{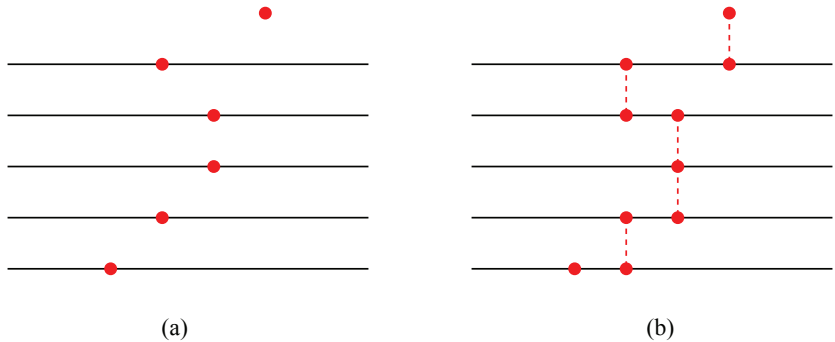}
\caption{The same function $\phi:=(\phi_i)_{i\geq 0}$ gives rise to two trap fields: (a) the model studied in~\cite{MOBC03, MOBC04}; (b) the model considered in Prop.~\ref{P:Pascal}.}
\label{fig:Trap}
\end{center}
\end{figure}

To prove Theorem~\ref{T:range}, it turns out to be fruitful to reformulate Theorem~\ref{T:range} in terms of a discrete time trapping problem. Note that analogous to the derivation of (\ref{rangecomp1}), we can rewrite
(\ref{range}) as
\be\label{target2}
\sum_{x\in\Z^d} \P^{\bar Z}_x( \tau_{-f}(\bar Z)\leq n) \geq \sum_{x\in\Z^d} \P^{\bar Z}_x( \tau_0(\bar Z)\leq n),
\ee
where
\be\label{tauf}
\tau_{-f}:=\tau_{-f}(\bar Z):=\min\{n\geq 0: \bar Z_n=-f_n\}.
\ee
We can reverse the roles of traps and particle and think of the trajectory $-f$ as a trap, with a particle starting from every site on $\Z^d$. The LHS of (\ref{target2}) is then the expected number of particles killed by the trap $-f$ by time $n$ (this is sometimes called a {\em target problem} in the physics literature, with $-f$ being the target). We can reformulate (\ref{target2}) in terms of a symmetric random walk $Z$ with increment distribution $p(\cdot)$ by contracting the time intervals $[2k, 2k+1]$, $k\geq 0$, into single time points. This results in a new trap field where at each time $i\geq 0$, there could be two sites in $\Z^d$ which act as traps (see Figure~\ref{fig:Trap}~(b)). More precisely, the new time-space trap field is given by
$$
\{ (i, -f_{2i}), (i, -f_{2i+1}) : i\geq 0\}.
$$
Denote $\phi_i=-f_{2i}$. Since $f_{2i+1}=f_{2i+2}$, the trap field equals $\{(i,\phi_i), (i,\phi_{i+1}) : i\geq 0\}$. For a path $Z:=(Z_n)_{n\geq 0}$ in $\Z^d$, denote
\be\label{ttauphi}
\tilde \tau_\phi:=\tilde \tau_\phi(Z) := \min\{n\geq 0: Z_n = \phi_n \mbox{ or } \phi_{n+1}\}.
\ee
Then we note that (\ref{target2}), with $n$ therein replaced by $2n+1$ for $n\in\N$, is equivalent to the following:
\bp\label{P:Pascal} {\bf [Pascal principle for discrete time trapping]}
Let $\P^Z_x$ denote probability for a random walk $Z$ on $\Z^d$ with $Z_0=x$, whose increment distribution $p(\cdot)$ satisfies one of the conditions in Theorem~\ref{T:range}~(i)--(iii). Then for any $\phi:=(\phi_i)_{i\geq 0}\in\Z^d$, we have
\be\label{target3}
\sum_{x\in\Z^d} \P^Z_x(\tilde\tau_\phi \leq n) \geq \sum_{x\in\Z^d} \P^Z_x(\tilde \tau_0\leq n) \qquad \mbox{for all } n\in\N.
\ee
\ep
When $n$ is even in (\ref{target2}), we can reduce it to the odd case of $n+1$ by setting $f_{n+1}:=f_n$. Thus we see that (\ref{target2}), and hence Theorem~\ref{T:range}, is equivalent to Proposition~\ref{P:Pascal}.

We now compare Proposition \ref{P:Pascal} with Moreau et al's result in~\cite{MOBC03, MOBC04}, which was briefly described in the introduction. The discrete time trapping problem considered by Moreau et al corresponds to a time-space trap field $\{(i, \phi_i): i\geq 0\}$ (see Figure~\ref{fig:Trap}~(a)), so that $\tilde\tau_\phi$ in (\ref{target3}) is replaced by $\tau_\phi$, defined as in (\ref{tauf}). The analogue of (\ref{target3}) is then proved in~\cite{MOBC03, MOBC04} for any random walk $Z$ whose increment distribution $p(\cdot)$ is symmetric with $p(0)\geq \frac{1}{2}$, and the corresponding monotonicity result for the random walk range is (see e.g.~\cite[Sec.~2.4]{DGRS12} for details)
\be\label{range2}
\E^Z_0[|R_n(Z-\phi)|] \geq \E^Z_0[|R_n(Z)|].
\ee
Note that the perturbation on the random walk path $Z$ is by adding a deterministic path $-\phi$ to $Z$, instead of inserting jumps in $Z$ as in Theorem~\ref{T:range}. Since $\tilde\tau_\phi \leq \tau_\phi$ and $\tilde\tau_0=\tau_0$, (\ref{target3}) with $\tilde\tau_\phi$ replaced by $\tau_\phi$ gives a stronger result, therefore Proposition~\ref{P:Pascal} and Theorem~\ref{T:range} hold for $p(\cdot)$ which are symmetric with $p(0)\geq \frac{1}{2}$.

However, as pointed out in~\cite{MOBC04}, (\ref{range2}) cannot hold for the simple symmetric random walk on $\Z^d$ because of periodicity. Indeed, if $d=1$ and we let $(\phi_n)_{n\geq 0}$ take values alternately between $0$ and $1$, then $\lim_{n\to\infty} \frac{|R_n(Z-\phi)|}{|R_n(Z)|}=\frac{1}{2}$ almost surely, because $Z-\phi$ can only visit even lattice sites. Therefore the insertion perturbation for the random walk range formulated in Theorem~\ref{T:range} is quite different from the additive perturbation considered in (\ref{range2}), and it is in a sense more natural since it includes the case of the simple symmetric random walk.

\section{A Class of Symmetric Random Walks on $\Z$}\label{S:Z}
In this section, we prove Proposition~\ref{P:Pascal} for symmetric random walks on $\Z$ with an increment distribution $p(\cdot)$ that
satisfies condition Theorem~\ref{T:range}~(i). Our initial goal was to prove Theorem~\ref{T:range} for the simple symmetric random walk on $\Z$, and in the process we discovered this more general class.

Our proof of Prop.~\ref{P:Pascal} for $p(\cdot)$, which satisfies condition Theorem~\ref{T:range}~(i), is based on induction and a suitable notion of symmetric domination of measures on $\Z$, which a priori may seem mysterious. A variant of this argument was recently applied by the second author with coauthors in~\cite{DSS11} to prove a rearrangement inequality for L\'evy processes in $\R^d$, which can be regarded as a generalized version of the Pascal principle for a trapping problem, where traps follow independent L\'evy motions. The argument in~\cite{DSS11} involves symmetric decreasing rearrangement of increment distributions in $\R^d$, which in the discrete setting we consider here would require $p(\cdot)$ to be symmetric and satisfy $p(0)\geq p(1)\geq p(2)\geq \cdots$. This however does not cover the simple symmetric random walk. Our key observation here is that by using a weaker notion of symmetric domination than in~\cite{DSS11} and careful manipulations, we can deal with $p(\cdot)$ that satisfies condition Theorem~\ref{T:range}~(i), which in particular includes the simple symmetric random walk.
\bigskip

\noindent
{\bf Proof of Prop.~\ref{P:Pascal} under condition Theorem~\ref{T:range}~(i).} For $n\geq 0$ and $x\in\Z$, define
\be\label{unx}
u^\phi_n(x) := 1-v^\phi_n(x) := \sum_{z\in\Z} \P^Z_z(Z_n=x, \tilde\tau_\phi>n),
\ee
and let $u^0_n(x)=1-v^0_n(x)$ denote the case $\phi_\cdot\equiv 0$. We can interpret $u^\phi_n(x)$ as the expected number of particles alive at time $n$ at position $x$, if initially one particle starts at every site in $\Z$ and gets killed when it encounters the time-space trap field $\{(i, \phi_i), (i, \phi_{i+1}) : i\geq 0\}$. Then (\ref{target3})
is equivalent to
\be\label{surviv}
\sum_{x\in\Z} v^0_n(x) \leq \sum_{x\in\Z} v^\phi_n(x).
\ee
We will prove (\ref{surviv}) by proving that $v^\phi_n$ symmetrically dominates $v^0_n$, denoted by $v^\phi_n\succ v^0_n$, in the sense that
\be\label{symdom}
\sum_{|x|\geq k} v^0_n(x) \leq \sum_{|x|\geq k} v^\phi_n(x_0+x) \qquad \forall \ k\geq 0,\ x_0\in\Z.
\ee
This notion of symmetric domination may seem mysterious. However, a heuristic interpretation is that, $u_n^0(\cdot)=1-v_n^0(\cdot)$ as a positive measure on $\Z$ puts more mass around $\infty$ and closer to $\infty$ than the measure $u_n^\phi(\cdot)=1-v_n^\phi(\cdot)$.
We will prove $v^\phi_n\succ v^0_n$, for all $n\geq 0$, by induction. It is easily verified that $v^\phi_0\succ v^0_0$. Now let us assume that $v^\phi_n\succ v^0_n$ for some $n\geq 0$ and try to prove that $v^\phi_{n+1} \succ v^0_{n+1}$.

Note that we have the following recursion relation for $u^\phi_n$:
\be
u^\phi_{n+1}(x) = \left\{
\begin{aligned}
0 \qquad \qquad & \qquad \qquad \  \mbox{if} \ x\in\{\phi_{n+1}, \phi_{n+2}\},\\
\sum_{y\in\Z} u^\phi_n(y) p(x-y) & \qquad \qquad \ \mbox{if}\ x\notin \{\phi_{n+1}, \phi_{n+2}\}.
\end{aligned}
\right.
\ee
Since $u^\phi_n(\phi_{n+1})=0$, we have
\be\label{vphi}
v^\phi_{n+1}(x) = 1-u^\phi_{n+1}(x) = \left\{
\begin{aligned}
1 \qquad \qquad \qquad & \qquad \  \mbox{if} \ x\in\{\phi_{n+1}, \phi_{n+2}\},\\
p(x-\phi_{n+1}) + \sum_{y\neq \phi_{n+1}} v^\phi_n(y) p(x-y) & \qquad \ \mbox{if}\ x\notin \{\phi_{n+1}, \phi_{n+2}\}.
\end{aligned}
\right.
\ee
Similarly,
\be\label{v0}
v^0_{n+1}(x) = \left\{
\begin{aligned}
1 \qquad \qquad \qquad & \qquad \  \mbox{if} \ x=0,\\
p(x) + \sum_{y\neq 0} v^0_n(y) p(x-y) & \qquad \ \mbox{if}\ x\neq 0.
\end{aligned}
\right.
\ee

We first verify the analogue of (\ref{symdom}) for $v^\phi_{n+1}$ and $v^0_{n+1}$, with $k=1$, which admits special cancelations and is the reason why we do not need to assume $p(0)\geq p(1)$. The case $k=0$ will be treated along the way.

Note that the second line in (\ref{vphi}) gives a lower bound on $v^\phi_{n+1}(x)$ for all $x\in\Z$, therefore
$$
\begin{aligned}
\sum_{|x|\geq 1} v^\phi_{n+1}(\phi_{n+1}+x) & \geq \sum_{|x|\geq 1}p(x) + \sum_{|x|\geq 1}\sum_{y\neq x} p(y) v^\phi_n(\phi_{n+1}+x-y) = F_1(0)+\sum_{z\neq 0} F_1(z) v^\phi_n(\phi_{n+1}+z), \\
\sum_{|x|\geq 1} v^0_{n+1}(x) & = \sum_{|x|\geq 1}p(x) + \sum_{|x|\geq 1}\sum_{y\neq x} p(y) v^0_n(x-y) = F_1(0)+\sum_{z\neq 0} F_1(z)v^0_n(z),
\end{aligned}
$$
where we made the change of variable $z=x-y$, and introduced the notation
$$
F_1(z):=\sum_{|x|\geq 1} p(z+x)= \sum_{|x|\geq 1} p(z-x) = \sum_{|y-z|\geq 1} p(y)=\sum_{|y+z|\geq 1} p(y).
$$
By the symmetry of $p(\cdot)$, $F_1(z)=F_1(-z)$, and by the assumption that $p(y)$ is decreasing in $y\geq 1$, we see that $F_1(z)$ is increasing in $z\geq 1$.
Indeed, for $z\geq 1$,
$$
F_1(z+1)-F_1(z) = \sum_{|y-z-1|\geq 1}p(y) - \sum_{|y-z|\geq 1}p(y) =p(z)-p(z+1)\geq 0.
$$
By a layer-cake representation for $F_1$, i.e., writing $F_1(z)=F_1(1)+\sum_{i=2}^{|z|} (F_1(i)-F_1(i-1))$, we have
\be\label{k1u0}
\begin{aligned}
& \sum_{|x|\geq 1} v^\phi_{n+1}(\phi_{n+1}+x) - \sum_{|x|\geq 1} v^0_{n+1}(x) \\
\geq \ \ & \sum_{z\neq 0} F_1(z)(v^\phi_n(\phi_{n+1}+z)-v^0_n(z)) \\
= \ \ & F_1(1)\sum_{|z|\geq 1} (v^\phi_n(\phi_{n+1}+z)-v^0_n(z)) + \sum_{i=2}^\infty(F_1(i)-F_1(i-1)) \sum_{|z|\geq i} (v^\phi_n(\phi_{n+1}+z)-v^0_n(z)),
\end{aligned}
\ee
which is non-negative because $v^\phi_n \succ v^0_n$ by the induction assumption.

Note that since $v^\phi_{n+1}(\phi_{n+1})=v^0_{n+1}(0)=1$, (\ref{k1u0}) implies
$$
\sum_x v^\phi_{n+1}(x)\geq \sum_x v^0_{n+1}(x),
$$
which in turn implies that
\be
\sum_{|x|\geq 1} v^\phi_{n+1}(x_0+x) \geq \sum_{|x|\geq 1} v^0_{n+1}(x) \qquad \forall\, x_0\in \Z,
\ee
because $v^\phi_{n+1}(x_0)\leq v^0_{n+1}(0)=1$. This verifies the analogue of (\ref{symdom}) for $v^\phi_{n+1}$ and $v^0_{n+1}$, with $k=0, 1$.

For $k\geq 2$ and $x_0\in \Z$, by (\ref{vphi}) and a change of variable $z:=x-y$, we have
$$
\sum_{|x|\geq k} v^\phi_{n+1}(x_0+x) \geq \sum_{|x|\geq k} \sum_{y\in\Z} p(y) v^\phi_n(x_0+x-y) = \sum_{z\in\Z} F_k(z) v^\phi_n(x_0+z),
$$
where
$$
F_k(z) := \sum_{|y+z|\geq k}p(y) = \sum_{|y-z|\geq k}p(y).
$$
Similarly, by (\ref{v0}) and the fact that $v^0_n(0)=1$, we have
$$
\sum_{|x|\geq k} v^0_{n+1}(x) = \sum_{z\in\Z} F_k(z) v^0_n(z).
$$
By the symmetry of $p(\cdot)$, $F_k(z)=F_k(-z)$, and for $z\geq 0$,
$$
F_k(z+1)-F_k(z) = \sum_{|y-z-1|\geq k}p(y)-\sum_{|y-z|\geq k}p(y) = p(z+1-k) - p(z+k).
$$
If $z+1-k\neq 0$, then $1\leq |z+1-k|< z+k$, and hence the above difference is non-negative by the assumption that $p(\cdot)$ is symmetric and $p(y)$ is decreasing in $y\geq 1$. If $z+1-k=0$, then the above difference is still non-negative because we are considering the case $k\geq 2$ and we assumed that $p(0)\geq p(3)$. Therefore $F_k$ is increasing in $z\geq 0$. As in (\ref{k1u0}), we can write
\be\label{ku}
\begin{aligned}
& \sum_{|x|\geq k} v^\phi_{n+1}(x_0+x) - \sum_{|x|\geq k} v^0_{n+1}(x) \\
\geq \ \ & \sum_{z\in\Z} F_k(z)(v^\phi_n(x_0+z)-v^0_n(z)) \\
= \ \ & F_k(0)\sum_{z\in\Z} (v^\phi_n(x_0+z)-v^0_n(z)) + \sum_{i=1}^\infty(F_k(i)-F_k(i-1)) \sum_{|z|\geq i} (v^\phi_n(x_0+z)-v^0_n(z)),
\end{aligned}
\ee
which again is non-negative because $v^\phi_n \succ v^0_n$ by the induction assumption. This completes the induction proof that $v^\phi_{n+1} \succ v^0_{n+1}$ defined as in (\ref{symdom}).
\qed

\section{The Simple Symmetric Random Walk on $\Z^d$}\label{S:SRW}

The argument we used in Section~\ref{S:Z} for random walks on $\Z$ unfortunately does not extend to higher dimensions. Instead, we
will follow a different line of argument to prove Proposition~\ref{P:Pascal} for the simple symmetric random walk in general dimensions.

Let $Z:=(Z_n)_{n\geq 0}$ be a random walk on $\Z^d$ with increment distribution $p(\cdot)$. Let $p_n(\cdot)$ denote the $n$-step increment distribution, i.e., $p_n(x) := \P^Z_0(Z_n=x)$, with $p_0(x)=\delta_0(x)$, and we set $p_{-1}(x):=0$ for all $x\in\Z^d$. For $n\geq -1$ and $x\in\Z^d$, we denote
$$
p_{n,n+1}(x):= p_n(x)+p_{n+1}(x).
$$
Prop.~\ref{P:Pascal} for the simple symmetric random walk on $\Z^d$ then follows from the following two lemmas.

\bl\label{L:dimd}
If $Z$ is a random walk whose transition probabilities satisfy
\begin{eqnarray}
p_{n,n+1}(0) &\geq& p_{n+1, n+2}(0) \qquad \forall \ n\geq -1, \label{mono1}\\
p_{n,n+1}(0) &\geq& p_{n, n+1}(x) \qquad \quad \forall \ n\geq -1, \, x\in\Z^d, \label{mono2}
\end{eqnarray}
then (\ref{target3}) holds.
\el

\bl\label{L:SRWmono}
If $Z$ is a simple symmetric random walk on $\Z^d$, then (\ref{mono1})--(\ref{mono2}) hold.
\el

\noindent
{\bf Proof of Lemma~\ref{L:dimd}.} Our proof of (\ref{target3}) is based on an adaptation of the argument of Moreau et al~\cite{MOBC03, MOBC04} for the discrete time trapping problem illustrated in Figure~\ref{fig:Trap}~(a). Although only an adaptation, to reach the conclusion that Moreau et al's argument can be adapted to our setting requires a crucial observation, which is non-trivial. Let us first recall the argument of Moreau et al, which is based on a simple but ingenious calculation.

To prove
\be\label{target4}
\sum_{x\in\Z^d} \P^Z_x(\tau_\phi \leq n) \geq \sum_{x\in\Z^d} \P^Z_x(\tau_0\leq n) \qquad \mbox{for all } n\in\N,
\ee
with $\tau_\phi$ and $\tau_0$ defined as in (\ref{tauf}), Moreau et al require the random walk to satisfy
\begin{eqnarray}
p_n(0) &\geq& p_{n+1}(0) \qquad \forall \ n\geq 0, \label{mono1'} \\
p_n(0) &\geq& p_n(x) \qquad  \quad \forall \ n\geq 0, \ x\in\Z^d, \label{mono2'}
\end{eqnarray}
which is easily seen to hold by Fourier transform if $p(\cdot)$ is symmetric and $p(0)\geq \frac{1}{2}$. Their idea is to consider a random walk $Z$ starting from any $x\in\Z^d$, restrict to the event that the walk ends at the trap location $\phi_n$ at time $n$, and perform a first passage decomposition according to the time when the walk first hits the trap $\phi_\cdot$. More precisely, they use the identity
\be\label{decomp1}
1=\sum_{x\in\Z^d} \P^Z_x(Z_n=\phi_n) = \sum_{i=0}^{n}  \sum_{x\in\Z^d} \P^Z_x(\tau_\phi =i)  p_{n-i}(\phi_n-\phi_i),
\ee
where the event $\{Z_n=\phi_n\}$ has been decomposed according to the value of $\tau_\phi$ (note that $Z_n=\phi_n$ implies $\tau_\phi\leq n$). Applying (\ref{mono2'}) in
(\ref{decomp1}), one obtains
$$
\sum_{i=0}^{n}  \sum_{x\in\Z^d} \P^Z_x(\tau_\phi =i)  p_{n-i}(0) \geq 1 = \sum_{i=0}^{n}  \sum_{x\in\Z^d} \P^Z_x(\tau_0 =i)  p_{n-i}(0).
$$
Define $W_\phi(-1):=0$, and $W_\phi(i):=W_\phi(i-1)+\sum_{x\in\Z^d} \P^Z_x(\tau_\phi=i)$ for $i\geq 0$. Define $W_0(\cdot)$ analogously. After rearranging terms, the above inequality can then be rewritten as
\be\label{tWcomp}
W_\phi(n)- W_0(n) \geq \sum_{i=0}^{n-1} \big(p_{n-i-1}(0)-p_{n-i}(0)\big)\,\big(W_\phi(i)- W_0(i)\big).
\ee
Since $p_{n-i-1}(0)-p_{n-i}(0)\geq 0$ for all $0\leq i\leq n-1$ by (\ref{mono1'}), one obtains $W_\phi(n)-W_0(n)\geq 0$ for all $n\geq 0$ by induction, which is precisely (\ref{target4}).

Note that the simple symmetric random walk does not satisfy (\ref{mono1'}) due to periodicity, and the above argument fails. Indeed, as explained at the end of Section~\ref{S:reduction}, (\ref{target4}) also fails for the simple symmetric random walk. Note that our model has two traps at each time, in contrast to Moreau et al's model, and Moreau et al's argument fails with no obvious remedy. We thus embarked on finding alternative approaches. However, after exhausting all other approaches, including the one we used in Section~\ref{S:Z} for dimension $1$, we returned to study Moreau et al's argument carefully once more, and we stumbled upon the following observation.

Although the simple symmetric random walk does not satisfy (\ref{mono1'}), it should satisfy (\ref{mono1}), while at the same time we note that in the trapping problem we are considering, each trap $\phi_i$ is present at both time $i$ and $i-1$ (the dashed lines in Figure~\ref{fig:Trap}~(b) serve to highlight this fact). Therefore instead of only considering the constraint $Z_n=\phi_n$ as in the decomposition (\ref{decomp1}), we can also consider the constraint $Z_{n-1}=\phi_n$. The event $Z_{n-1}=\phi_n$ implies $\tilde\tau_\phi \leq n-1$, which allows us to obtain decompositions analogous to (\ref{decomp1}):
\be\label{decomp2}
\begin{aligned}
1 & = & \sum_{x\in\Z^d} \P^Z_x(Z_n=\phi_n) & = \sum_{i=0}^{n}  \sum_{x\in\Z^d} \E^Z_x\big[1_{\{\tilde\tau_\phi =i\}} p_{n-i}(\phi_n-Z_i)\big], \\
1 & = & \sum_{x\in\Z^d} \P^Z_x(Z_{n-1}=\phi_n) & = \sum_{i=0}^{n-1}  \sum_{x\in\Z^d} \E^Z_x\big[1_{\{\tilde\tau_\phi =i\}} p_{n-1-i}(\phi_n-Z_i)\big].
\end{aligned}
\ee
Note that given $\tilde\tau_\phi=i$, $Z_i\in \{\phi_i, \phi_{i+1}\}$. Recall that $p_{-1}(\cdot):= 0$, adding the two equalities then gives
$$
2 = \sum_{i=0}^{n}  \sum_{x\in\Z^d} \E^Z_x\big[1_{\{\tilde\tau_\phi =i\}} \big(p_{n-1-i}(\phi_n-Z_i)+p_{n-i}(\phi_n-Z_i)\big)\big].
$$
This identity holds in particular for the case $\phi\equiv 0$. We can now apply condition (\ref{mono2}) to obtain
\be
\sum_{i=0}^{n}  \sum_{x\in\Z^d} p_{n-i-1, n-i}(0) \P^Z_x(\tilde\tau_\phi =i) \geq 2 = \sum_{i=0}^{n}  \sum_{x\in\Z^d} p_{n-i-1, n-i}(0) \P^Z_x(\tilde\tau_0 =i).
\ee
Define $\tilde W_\phi(-1):=0$, and $\tilde W_\phi(i):=\tilde W_\phi(i-1)+\sum_{x\in\Z^d} \P^Z_x(\tilde \tau_\phi=i)$ for $i\geq 0$. Define $\tilde W_0(\cdot)$ analogously. Rearranging terms in the above inequality then gives the following analogue of (\ref{tWcomp}):
\be\label{tWcomp2}
\tilde W_\phi(n)- \tilde W_0(n) \geq \sum_{i=0}^{n-1} \big(p_{n-i-2, n-i-1}(0)-p_{n-i-1,n-i}(0)\big)\,\big(\tilde W_\phi(i)- \tilde W_0(i)\big).
\ee
By assumption (\ref{mono1}), $p_{n-i-2,n-i-1}(0)-p_{n-i-1,n-i}(0)\geq 0$ for all $0\leq i\leq n-1$, which implies that $\tilde W_\phi(n)-\tilde W_0(n)\geq 0$ for all $n\geq 0$ by induction. This is precisely (\ref{target3}).
\qed
\bigskip

\noindent
{\bf Proof of Lemma~\ref{L:SRWmono}.} For $k=(k_1,\cdots, k_d)\in [-\pi,\pi]^d$, let $\psi(k):=\E^Z_0[e^{i k\cdot Z_1}]=\frac{1}{d}\sum_{i=1}^d \cos k_i$. Then
\be\label{pnx}
p_n(x) = \frac{1}{(2\pi)^d} \int_{[-\pi, \pi]^d} e^{-i k\cdot x} \psi^n(k) {\rm d}k \qquad \forall\ x\in\Z^d.
\ee
By periodicity, $p_n(0)=0$ for all $n$ odd, while the above identity shows that $p_{2n}(0)$ is decreasing in $n$. These facts readily imply (\ref{mono1}).

Using (\ref{pnx}), we can easily deduce (\ref{mono2}) for the case $n$ is even. However when $n$ is odd, (\ref{mono2}) does not seem to admit a simple proof using characteristic functions. Instead we will proceed via a coupling argument. Assume that $n\in\N$ is odd. If $x=(x_1,\cdots, x_d)\in\Z^d$ is such that $|x|_1:=\sum_{i=1}^d |x_i|$ is even, then $p_{n,n+1}(x)=p_{n+1}(x)\leq p_{n+1}(0)=p_{n,n+1}(0)$, where the inequality can be deduced from (\ref{pnx}). If $|x|_1$ is odd, then $p_{n,n+1}(x)=p_n(x)$. Also note that $p_{n,n+1}(0)=p_{n+1}(0)=p_n(e_1)$, where $e_1=(1,0,\cdots,0)\in\Z^d$. Therefore to complete the proof of (\ref{mono2}),
it only remains to show that
\be\label{pnxodd}
p_n(x) \leq p_n(e_1) \qquad \forall\ n\in \N, \ x\in\Z^d \mbox{ with } n \mbox{ and } |x|_1 \mbox{ both odd}.
\ee
By the symmetry of the random walk, $p_n(x)=\P^X_x(X_n=0)$ and $p_n(e_1) = \P^Y_{e_1}(Y_n=0)$ for two simple symmetric random walks $X$ and $Y$, starting respectively at
$x$ and $e_1$. We will construct a coupling of $X$ and $Y$ such that when $X_n=0$, we also have $Y_n=0$.

By symmetry, we may assume without loss of generality that $x_i\geq 0$ for all $1\leq i\leq d$. Since $|x|_1$ is odd, we may further assume that $x_i$
is odd for $1\leq i\leq 2m-1$ for some $m\in\N$, and $x_i$ is even for $2m\leq i\leq d$. We will group $(x_i)_{2\leq i\leq 2m-1}$ into $m-1$ pairs: $(x_2, x_3),\cdots, (x_{2m-2}, x_{2m-1})$. Let $X_k^{(i)}$ and $Y_k^{(i)}$ denote respectively the $i$-th component of $X_k$ and $Y_k\in\Z^d$. We will couple $X=(X_k^{(1)},\cdots, X_k^{(d)})_{k\geq 0}$ and $Y=(Y_k^{(1)},\cdots, Y_k^{(d)})_{k\geq 0}$ in such a way that:
\begin{itemize}
\item[(1)] at each time $k\geq 0$, $|X_k^{(i)}|\geq |Y_k^{(i)}|$ for all $1\leq i\leq d$;
\item[(2)] if $X_k^{(i)}=Y_k^{(i)}$, then $X_{k'}^{(i)}=Y_{k'}^{(i)}$ for all $k'\geq k$;
\item[(3)] if $X_k^{(i)}\equiv Y_k^{(i)} \mbox{ (mod 2)}$, then $X_{k'}^{(i)}\equiv Y_{k'}^{(i)} \mbox{ (mod 2)}$ for all $k'\geq k$;
\item[(4)] for each $1\leq i\leq m-1$, $X_k^{(2i)}+Y_k^{(2i)}  \equiv X_k^{(2i+1)} +Y_k^{(2i+1)} \mbox{ (mod 2)}$ for all $k\geq 0$.
\end{itemize}
Such a coupling clearly would imply $\P^X_x(X_n=0)\leq \P^Y_{e_1}(Y_n=0)$. In words, our coupling will be such that: if $X_k^{(i)}=Y_k^{(i)}$ for some $1\leq i\leq d$, then we couple the jumps of $X$ and $Y$ in the $i$-th coordinate so that $X^{(i)}_\cdot=Y^{(i)}_\cdot$ for all later times; if $X_k^{(i)}\neq Y_k^{(i)}$ but they have the same parity, then we couple the jumps in the $i$-th coordinate so that $X^{(i)}_\cdot$ and $Y^{(i)}_\cdot$ move as mirror images across $\frac{X_k^{(i)}+Y_k^{(i)}}{2}$, until the first time $X^{(i)}_\cdot=Y^{(i)}_\cdot$; if $X_k^{(i)}$
and $Y_k^{(i)}$ do not have the same parity, then $2\leq i\leq 2m-1$, and by (4), we can find another coordinate $i'$ such that $X_k^{(i')}$ and $Y_k^{(i')}$ also have different parities, and we couple the jumps in the $i$-th and $i'$-th coordinates such that the jump of $X$ in the $i$-th coordinate coincides with the jump of $Y$ in the $i'$-th coordinate and vice versa. The precise formulation is as follows.

Assume that $X$ and $Y$ have been coupled up to time $k$, and properties (1)--(4) have not been violated. Then $X_k^{(i)}\equiv Y_k^{(i)} \mbox{ (mod 2)}$ for $i=1$ and $2m\leq i\leq d$. Given $X_{k+1}-X_k\in \{e_i, -e_i\}$ for some $1\leq i\leq d$:
\begin{itemize}
\item[(a)] if $X_k^{(i)}=Y_k^{(i)}$, then we set $Y_{k+1}:=Y_k+ (X_{k+1}-X_k)$;
\item[(b)] if $X_k^{(i)}\equiv Y_k^{(i)} \mbox{ (mod 2)}$ and $X_k^{(i)}\neq Y_k^{(i)}$, then we set $Y_{k+1}:=Y_k - (X_{k+1}-X_k)$;
\item[(c)] if $X_k^{(i)}\not\equiv Y_k^{(i)} \mbox{ (mod 2)}$, then $i\in \{2j, 2j+1\}$ for some $1\leq j\leq m-1$. Let $i'\in\{2j, 2j+1\}\backslash\{i\}$.
Then we set $Y_{k+1}:=Y_k + e_{i'} \langle e_i, X_{k+1}-X_k\rangle$ where $\langle,\rangle$ denotes inner product on $\R^d$.
\end{itemize}
Note that $Y$ is distributed as a simple symmetric random walk. Conditions (2)--(4) remain un-violated at time $k+1$. We now check that (1) also holds at time $k+1$.
Indeed, jumps as specified in (a) do not violate (1). For a jump as specified in (b), by our assumptions $|X_k^{(i)}|\geq |Y_k^{(i)}|$ and $X_k^{(i)}\equiv Y_k^{(i)} \mbox{ (mod 2)}$, either $X_k^{(i)}=-Y_k^{(i)}$, in which case $X_{k+1}^{(i)}=-Y_{k+1}^{(i)}$; or $|X_k^{(i)}|\geq |Y_k^{(i)}|+2$, in which case we must have $|X_{k+1}^{(i)}|\geq |Y_{k+1}^{(i)}|$. For a jump as specified in (c), we must have $|X_k^{(i)}|\geq 1+|Y_k^{(i)}|$, and by (4), also $|X_k^{(i')}|\geq 1+|Y_k^{(i')}|$. Since $X$ and $Y$ do not make jumps in the same coordinate,
(1) must hold at time $k+1$ as well. This verifies that the coupling given in (a)--(c) satisfies conditions (1)--(4) at all times, which concludes the proof of the lemma.
\qed

\bigskip

\noindent
{\bf Acknowledgement} We thank Alex Drewitz for helpful comments on an earlier version of this paper. L.-C.~Chen is supported by research grant 99-2115-M-030-004-MY3 from the National Science Council of Taiwan, and he thanks the National University of Singapore for support during research visits. R.~Sun is supported by research grant R-146-000-119-133 from the National University of Singapore, and he thanks Academia Sinica of Taiwan for support during research visits.

\end{document}